\documentclass[11pt]{amsart}

\usepackage{color}
\usepackage{amsmath}
\usepackage{amssymb}
\usepackage{latexsym}
\usepackage{amsfonts}

\def\a{\alpha}
\def\b{\beta}
\def\d{\delta}

\def\g{\gamma}

\def\t{\tau}

\def\om{\Omega}

\def\ts{\times}
\def\iy{\infty}

\def\im{{\rm Im\ }}

\def\ker{{\rm Ker\,}}
\def\diag{{\rm diag\,}}

\def\BC{{\mathbb C}}

\def\BR{{\mathbb R}}

\def\cla{{\mathcal A}}
\def\clb{{\mathcal B}}
\def\clc{{\mathcal C}}
\def\cld{{\mathcal D}}

\def\clp{{\mathcal P}}

\newcommand{\kr}{{\rm Ker\,}}

\def\bfr{\mathbf{R}}
\def\bft{\mathbf{T}}

\newcommand{\bpr}{{\noindent\textbf{Proof.}\ \ }}
\newcommand{\epr}{{\hfill $\Box$}}

\newcommand{\la}{\lambda}

\newtheorem{Pa}{Paper}[section]
\newtheorem{Tm}[Pa]{{\bf Theorem}}
\newtheorem{La}[Pa]{{\bf Lemma}}
\newtheorem{Cy}[Pa]{{\bf Corollary}}

\newtheorem{Pn}[Pa]{{\bf Proposition}}

\title{Krein systems}

\author[Alpay]{D. Alpay}
\address{(DA) Department of Mathematics\\
Ben--Gurion University of the Negev\\
Beer-Sheva 84105\\ Israel}
\email{dany@math.bgu.ac.il}

\author[Gohberg]{I. Gohberg}
\address{(IG)
School of Mathematical Sciences\\
The Raymond and Beverly Sackler
Faculty of Exact Sciences\\
Tel--Aviv University \\
Tel--Aviv, Ramat--Aviv 69989, Israel}
\email{gohberg@post.tau.ac.il}

\author[Kaashoek]{M.A. Kaashoek}
\address{(MK)
Afdeling Wiskunde\\
Faculteit der Exacte Wetenschappen\\
Vrije Universiteit\\
De Boelelaan 1081a, 1081 HV Amsterdam, The
Netherlands}
\email{ma.kaashoek@few.vu.nl}

\author[Lerer]{L. Lerer}
\address{(LL)
Department of Mathematics\\
Technion, Israel Institute of Technology\\
Haifa 32000, Israel}
\email{llerer@techunix.technion.ac.il}

\author[Sakhnovich]{A. Sakhnovich}
\address{(AS)
Fakult\"at f\"ur Mathematik\\
Universit\"at Wien\\
Nordbergstrasse 15, A-1090 Wien, Austria}
\email{al$\_$sakhnov@yahoo.com}

\thanks{Daniel Alpay wishes to thank the
Earl Katz family for endowing the chair which
supported his research. The work of Alexander Sakhnovich was supported by the
Austrian Science Fund (FWF) under
Grant  no. Y330}

\dedicatory{In memory of Mark Grigorievich Krein, with appreciation\\ of his
many great discoveries, on the occasion of his Centennial.}

\begin{document}

\date{}

\begin{abstract}
In the present paper we extend
results of M.G. Krein associated to the spectral
problem for Krein systems to systems with matrix valued accelerants
with a possible jump discontinuity at the origin. Explicit formulas
for the accelerant are given in terms of
the matrizant of the system in question. Recent developments in
the theory of continuous analogs of the
resultant operator play an essential role.
\end{abstract}
\maketitle

\section{Introduction}
\setcounter{equation}{0}

\noindent
The following result is due to M.G. Krein, see
\cite{MR0086978}:
\begin{Tm}
\label{240308} Let ${\mathbf T}>0$, and let $k$
be a scalar continuous and hermitian function on
the interval $[-\mathbf T,\mathbf T]$ such that
for each $0< \t\leq \bft$ the corresponding
convolution integral operator $T_\t$ on
${\mathbf L}^1[0,\,\t]$,
\begin{equation}
\label{intor} (T_\t f)(t)=f(t)-\int_0^\t
k(t-s)f(s)\,ds, \quad 0\leq t\leq \t,
\end{equation}
is invertible. Let $\gamma_\t(t,s)$ denote the
resolvent kernel
\begin{equation}
\label{langoustine}
\gamma_\tau(t,s)- \int_0^\tau
k(t-v)\gamma_\tau(v,s)dv=k(t-s),\quad 0\le t,s\le
\tau.
%
%
\end{equation}
Consider the entire function
\begin{eqnarray}
\label{defP} {\mathcal P}(\tau,\la)&=& e^{i\la
\tau} \left(1+\int_0^{\tau}e^{-i\la
x}\gamma_{\tau}(x,0)
dx\right),\\
\label{defP*} {\mathcal
P}_*(\tau,\la)&=&1+\int_0^{\tau}e^{i\la
x}\gamma_{\tau}(\tau-x, \tau) dx.
\end{eqnarray}
Then with $a(\t)=\g_\t(\t,0)$ and for $\la\in\BC$ it holds that
\begin{equation}
\label{kr1}
\left\{
\begin{array}{l}
\displaystyle{\frac{\partial}{\partial\t}}
\clp(\t,\la)= i\la \clp(\t,\la)
+\clp_{*}(\t,\la)a(\t),
\quad 0\leq \t\leq\bft,\\
\noalign{\vskip6pt}
\displaystyle{\frac{\partial}{\partial\t}}
\clp_*(\t,\la)= \clp(\t,\la)a(\t)^*.
\end{array}
\right.
\end{equation}
\end{Tm}

\medskip

Putting $Y(\t,\la)=\begin{bmatrix} {\mathcal
P}(\t,\la)&
{\mathcal P}_*(\t,\la)\end{bmatrix}$,  the system \eqref{kr1}
can be rewritten as
\begin{equation}
\label{kr} \frac{\partial}{\partial \t}
Y(\t,\la)=Y(\t,\la)
\left(i\la\begin{bmatrix}I_r&0\\0&0
\end{bmatrix}+\begin{bmatrix}0&a(\t)\\ a(\t)^*&
0\end{bmatrix}\right).
\end{equation}
Here  $\t\in[0,\bft]$. We call (\ref{kr}) a \emph{Krein system}
when, as in (\ref{kr1}), the function  $a$ is given by  $a(\t)=\gamma_\t(\t,0)$,
where $\g_\t(t,s)$ is the resolvent kernel corresponding to
some $k$ on $[-\bft, \bft]$ with the properties described in the Theorem \ref{240308}.
In that case, following Krein,  the function $k$ is called an \emph{accelerant} for  (\ref{kr}),
and we shall refer to $a$  as the \emph{potential associated with the
accelerant} $k$. The functions $\clp(\t,\cdot)$, $\clp_*(\t,\cdot)$ are  called
\emph{Krein orthogonal functions} at $\t$ associated to the \emph{weight}  $\d-k$, where $\d$ is the delta function.

In this paper we prove the analogue of Theorem \ref{240308} for
systems with
accelerants that are allowed to have a jump
discontinuity at the origin. We also present explicit formulas
for determining the unique accelerant $k$ from
the given potential $a$. As for  continuous accelerants in \cite{MR0086978},
the results are proved not only for  scalar
functions but also
for the matrix-valued case,  when in \eqref{kr1} the functions
${\mathcal P}$, ${\mathcal P}_*$ and $a$ are
${\mathbb C}^{r\times r}$--valued.

The result expressing the accelerant in terms of the potential   referred to
in the previous paragraph is based on a recent theorem involving
a certain analog $\bfr(B,D)$ of the resultant operator for a
class of entire matrix functions $B$ and $D$.
The resultant $\bfr(B,D)$ is defined as follows (see
Section \ref{270308} for more details). Let $B$
and $D$ be of the form
\[
B(\la)=I_r+\int_{-\t}^0e^{i\la u}b(u)du\quad{\rm
and}\quad D(\la)=I_r+\int^{\t}_0e^{i\la u}d(u)du,
\]
where the functions $b$ and $d$ belong
respectively to ${\mathbf L}_1^{r\times
r}[-\t,0]$ and ${\mathbf L}_1^{r\times r}[0,\t]$.
The {\sl resultant} of $B$ and $D$ is the
operator defined on the space ${\mathbf L}_1^{r\times
r}[-\t,\t]$ by:
\[
(\bfr(B,D)q)(u)=\begin{cases}q(u)+ \displaystyle{\int_{-\t}^\t
d(u-s)q(u)du},\quad 0\le
u\le \t,\\
\noalign{\vskip6pt}
q(u)+\displaystyle{\int_{-\t}^\t b(u-s)q(u)du},\quad -\t\le u<
0.
\end{cases}
\]
Let us now state our  main results.

\begin{Tm}
\label{mainthm1} Let $k$ be  a ${r\ts r}$-matrix valued accelerant on $[-\mathbf
T,\mathbf T]$, with possibly a jump
discontinuity at the origin,  and let $\g_\t(t,s)$ be the corresponding resolvent kernel as
in $(\ref{langoustine})$. Put
\begin{eqnarray}
\label{defP_mat} {\mathcal P}(\tau,\la)&=&
e^{i\la \tau} \left(I_r+\int_0^{\tau}e^{-i\la
x}\gamma_{\tau}(x,0)
dx\right)\\
\label{defP*_mat} {\mathcal
P}_*(\tau,\la)&=&I_r+\int_0^{\tau}e^{i\la
x}\gamma_{\tau}(\tau-x, \tau) dx.
\end{eqnarray}
Then $a(\tau)=\gamma_\tau(0,\tau)$,
with $0<\tau\le \mathbf T$, extends to a continuous function on $[0,\bft]$ and
\[
Y(\t,\la)=\begin{bmatrix} {\mathcal
P}(\t,\la)&
{\mathcal P}_*(\t,\la)\end{bmatrix},
\]
satisfies the Krein system \eqref{kr} with
potential $a$.
\end{Tm}

For our second main result we need the \emph{matrizant} of \eqref{kr}.
By definition, this is the unique ${\mathbb
C}^{2r\times 2r}$--valued solution $U(\t,\la)$ of \eqref{kr} satisfying the initial condition
$U(0,\la)\equiv I_{2r}$.

\begin{Tm}
\label{mainthm2}Let $k$ be  a ${r\ts r}$-matrix valued accelerant on $[-\mathbf
T,\mathbf T]$, with possibly a jump
discontinuity at the origin, and let $a$ be the corresponding
potential. Then $k$ is uniquely determined by $a$, and $k$ can be obtained from $a$
in the following way. Let $U(\t,\la)$ be the matrizant of \eqref{kr}, and put
\[
F(\la)=e^{i\la\bft}\begin{bmatrix} I_r&I_r
\end{bmatrix}U(\bft, -\la)\begin{bmatrix}
I_r\\0
\end{bmatrix}, \  G(\la)=\begin{bmatrix} I_r&I_r
\end{bmatrix}U(\bft, -\la)\begin{bmatrix}
0\\I_r
\end{bmatrix}.
\]
Then $F$ and $G$ are entire $r\ts r$ matrix functions of the form
\[
F(\la)= I_r+\int_0^{\bft}f(x)e^{i\la
x}\,dx, \quad G(\la)=I_r+\int_{-\bft}^0 g(x)e^{i\la
x}\,dx,
\]
where $f$ and $g$ are continuous $\BC^{r\ts r}$-valued functions on
$[0,\bft]$ and $[-\bft, 0]$, respectively. Moreover, the resultant operator
$\bfr(F^\sharp,G^\sharp)$ is invertible, and
the function $k$ is given by the formula
\begin{equation}
\label{Pernety}
k=[\bfr(F^\sharp,G^\sharp)]^{-1}q.
\end{equation}
Here  $F^\sharp(\la)=F(\bar{\la})^*$ and $G^\sharp(\la)=G(\bar{\la})^*$, where the superscript ${}^*$
means taking adjoints. Finally, $q$  is the function on the interval $[-\bft, \bft]$ given by
\[
q(x)=\begin{cases}
f(-x)^*,\quad -\bft\le x<0,\\
g(-x)^*,\quad \hspace{4mm}0\le x\le \bft.
\end{cases}
\]
\end{Tm}

To prove Theorem \ref{mainthm1} we use in an
essential way the results of \cite{goh-kol}.
The proof  of Theorem \ref{mainthm2} is based on
recent
results of \cite{GKLe07b} on the continuous
analog of the resultant.

In each of the two  theorems above our starting point is a given accelerant.
In a next paper we plan to study the inverse situation, which includes, in particular, the
question
whether or not a continuous potential is always generated by
an accelerant.

Let us illustrate Theorem \ref{mainthm1} with an example. Take $k$ to be
\begin{equation}
\label{Ranelagh}
k(t)=\begin{cases}
\,\,\,\,i,\quad{\rm if}\quad t\in[0,{\mathbf T}],\\
-i,\quad{\rm if}\quad t\in[-{\mathbf T},0].
\end{cases}
\end{equation}
Clearly, $k$ is continuous with a jump discontinuity at zero, and  $k$ is hermitian.
Note that this function
is of the form
\begin{equation}
\label{expkernelintro}
k(t)=\begin{cases}iCe^{-itA}
(I-P)B,\quad t\in[0,{\mathbf T}],\\
-iCe^{-itA}PB,\,\,\,\,\qquad t\in[-{\mathbf
T},0],
\end{cases}
\end{equation}
with
\[
A=\begin{bmatrix}0&0\\0&0\end{bmatrix},\quad
P=\begin{bmatrix}0&0\\0&1\end{bmatrix},\quad
{\rm and}\quad
C=B^*=\begin{bmatrix}1&1\end{bmatrix}.
\]
The formulas from \cite{bgk3} allow us to show that for this $k$
the integral operator $T_\t$ in \eqref{intor} is invertible  for ${\t}<\frac{\pi}{2}$. Hence
$k$  is an accelerant on $[-\bft,\bft]$ whenever ${\mathbf T}<\frac{\pi}{2}$.
Furthermore, again using the formulas from \cite{bgk3}, one computes
that for each ${\t}<\frac{\pi}{2}$  the resolvent kernel associated to $k$, that is,
the solution $\gamma_\t(t,s)$ of
\eqref{langoustine}, is given by
\[
\gamma_\t(t,s)=\begin{cases}
\dfrac{ie^{2i(t-s)}}{1+e^{2i\t}},\quad
0\leq s<t\le{\t},
\\
\mbox{}
\\
\dfrac{-ie^{2i(t-s)}}{1+e^{-2i\t}},\quad
0\leq t<s\le{\t}.
\end{cases}
\]
Direct
computations show then that the functions
${\mathcal P}$ and ${\mathcal P}_*$ defined by
the formulas \eqref{kr1} are equal to
\[
\begin{split}
{\mathcal
P}(\t,\la)&=e^{i\la\t}+\dfrac{2}{1+e^{2i\t}}\dfrac{
e^{2i\t}-e^{i\la\t}}{2-\la},\\
{\mathcal
P}_*(\t,\la)&=1+\dfrac{2}{1+e^{2i\t}}\dfrac{
e^{2i\t}-e^{i\la\t}}{2-\la},
\end{split}
\]
and that these functions satisfy the system
\eqref{kr1} with
\begin{equation}
\label{Pont_Marie}
a(\tau)=\dfrac{2i}{1+e^{-2i\t}},
\quad\tau\in[0,{\mathbf T}].
\end{equation}
Other examples will be given in the final two sections of the paper.

We now give the outline of the paper. The rest
of the paper consists of    five sections. In
Section \ref{vienne} we show that a Krein system can be
associated to accelerants with jump
discontinuities and prove Theorem \ref{mainthm1}. In Section \ref{270308} we review
the notion of continuous analogue of the
resultant and state the results from \cite{GKLe07b} used in this paper. The proof of Theorem
\ref{mainthm2} is given in Section \ref{seckrs}. The last two
sections present examples. In Section \ref{exjump} we
consider the case of accelerants $k$ of the form \eqref{expkernelintro},
where $A,B$ and $C$ are matrices
of appropriate sizes and $P$ is a projection
commuting with $A$. This includes in
particular the case when the Fourier
transform of $k$ (considered as a function on $\BR$)
is a rational matrix--valued
function vanishing at infinity. Such functions
$k$ have in general a jump discontinuity at the
origin. In Section \ref{otherclass} a class of continuous
accelerants
is elaborated.

\section{Krein system for accelerants with
jump discontinuity\\ and proof of Theorem \ref{mainthm1}}
\label{vienne}
\setcounter{equation}{0}
\noindent
In the proof of Theorem \ref{240308} an important role is played
by the equations
\begin{eqnarray} \frac{\partial}{\partial\t}\,
\g_\t(t,s)&=&\g_\t(t,\t)\g_\t(\t,s), \quad 0\leq
s, t\leq \t,\label{KSeq1}\\
\frac{\partial}{\partial\t}\,
\g_\t(\t-t,\t-s)&=&\g_\t(\t-t,0)\g_\t(0,\t-s),
\quad 0\leq t, s \leq \t.\label{KSeq2}
\end{eqnarray}
Equation \eqref{KSeq1} is called the
\emph{Krein-Sobolev identity}. The second equation is
obtained by replacing in equation
\eqref{langoustine} the function $k(t)$ by $k(-t)$. The
corresponding resolvent kernel is equal to
\(\gamma_\tau(\tau-t,\tau-s)\), as can be seen
by a change of variables; see the discussion
\cite[p. 450]{dym-90} and in particular equation
$(3.5)$ there. The above equations have been
used by M.G. Krein in \cite{Krein9} to deduce his
system \eqref{kr1} in the case of a continuous
accelerant.

It is known  \cite[Section 7.3, p.\,187]{gk-2} that
continuity of the accelerant is not necessary to
insure that the Krein-Sobolev identity holds. In fact, when $k$
 has a jump discontinuity at the origin appropriate
generalizations of \eqref{KSeq1}-\eqref{KSeq2}
have been established in
\cite{goh-kol}.

Before presenting the proof of Theorem \ref{mainthm1}, we
first review the necessary results from \cite{goh-kol}. In what follows $k$
is a $r\ts r$ accelerant on $[-\bft,\bft]$ with a possible jump
discontinuity at the origin and $\g_\t(t,s)$ is the corresponding
resolvent kernel as in \eqref{langoustine}. From \cite{goh-kol} we
know that the function
$(t,s,\t)\mapsto \gamma_\t(t,s)$ is continuous on the domain $0\leq s<t\leq \bft,\ 0< \t\leq \bft$
and on the domain $0\leq t<s\leq \bft,\ 0< \t\leq \bft$. Moreover,
$(t,s,\t)\mapsto \gamma_\t(t,s)$ admits continuous extensions on the closures of these domains.
In particular, $a(\t)=\g_\t(\t,0)$ is continuous on the left open interval $(0,\bft]$ and
has a continuous extension to
the closed interval $[0,\bft]$.

Next, we consider the modifications of equations \eqref{KSeq1}--\eqref{KSeq2}.
Using the fact that $k$ has a jump discontinuity at the origin, we let $k_+$ be the
function equal to $k$ for $t\not =0$
\[
k_+(0)=\lim_{\substack{ h\rightarrow 0\\
h>0}}k(h).
\]
Similarly, let $k_-$ be the function equal to $k$
for $t\not =0$ and defined at the origin by
\[
k_-(0)=\lim_{\substack{ h\rightarrow 0\\
h<0}}k(h).
\]
One defines $\gamma_\tau^{u}(t,s)$ and
$\gamma_\tau^l(t,s)$ to be the resolvent
equations corresponding to the function $k_+(t)$
and $k_-(t)$ respectively. Note that, for
$t\not=s$,
\begin{equation}
\label{paris}
\gamma_\tau(t,s)=\gamma^u_\tau(t,s)=
\gamma^l_\tau(t,s).
\end{equation}
%
It is proved in \cite{goh-kol} that
\begin{eqnarray}
\frac{\partial}{\partial\t^+}\,
\g_\t^u(t,s)&=&\g_\t^u(t,\t)\g_\t^l(\t,s),
\quad 0\leq s, t\leq \t,\label{KSeq11}\\
\frac{\partial}{\partial\t^+}\,
\g_\t^l(t,s)&=&\g_\t^u(t,\t)\g_\t^l(\t,s), \quad
0\leq s, t\leq \t,\label{KSeq111}
\end{eqnarray}
and
\begin{eqnarray}
\frac{\partial}{\partial\t^-}\,
\g_\t^u(t,s)&=&\g_\t^u(t,\t)\g_\t^l(\t,s),
\quad 0\leq s, t\leq \t,\label{KSeq117}\\
\frac{\partial}{\partial\t^-}
\g_\t^ls(t,s)&=&\g_\t^u(t,\t)\g_\t^l(\t,s)\,\quad
0\leq s, t\leq \t,\label{KSeq1117}
\end{eqnarray}
where $\frac{\partial}{\partial^+}$ and
$\frac{\partial}{\partial^-}$ stand for
derivatives from the right and from the left,
respectively. See \cite[(3.6)-(3.7) p.\,274, and
p.\,278]{goh-kol}. It follows that \eqref{KSeq2}
becomes
\begin{eqnarray}
\frac{\partial}{\partial\t^+}\,\g_\t^u(\t-t,\t-s)&=
&\g^u_\t(\t-t,0)\g^l_\t(0,\t-s), \label{KSeq22}
\\
\frac{\partial}{\partial\t^+}\,\g_\t^l(\t-t,\t-s)&=
&\g^u_\t(\t-t,0)\g^l_\t(0,\t-s), \label{KSeq222}
\end{eqnarray}
where $0\leq t,s \leq \t$, and
\begin{eqnarray}
\frac{\partial}{\partial\t^-}\,\g_\t^u(\t-t,\t-s)&=
&\g^u_\t(\t-t,0)\g^l_\t(0,\t-s), \label{KSeq2267}
\\
\frac{\partial}{\partial\t^-}\,\g_\t^l(\t-t,\t-s)&=
&\g^u_\t(\t-t,0)\g^l_\t(0,\t-s),
\label{KSeq22222}
\end{eqnarray}
also for $0\leq t,s \leq \t$.\\

\noindent\textbf{Proof of Theorem \ref{mainthm1}.} We have already proved
the continuity of the potential $a$ on $[0,\bft]$.

Let $\clp$ and $\clp_*$ be defined by \eqref{defP_mat} and \eqref{defP*_mat}.
Note that, in view of \eqref{paris},
one can replace $\gamma_\tau$ by $\gamma_\tau^u$
or $\gamma_\tau^l$ in the expressions for $\clp$
and $\clp_*$. Then, using the Krein-Sobolev
identity (\ref{KSeq11}), we have for $\tau>0$:
\begin{eqnarray*}
\frac{\partial}{\partial\t^+}\clp(\t,\la)&=& i\la
\clp(\t,\la)+e^{i\la\t}
\frac{\partial}{\partial\t^+}\int_0^\t
e^{-i\la x}\g_\t(x,0)\,dx\\
&=& i\la \clp(\t,\la)+\g_\t(\t,0)+ \int_0^\t
e^{i\la(\t- x)}
\frac{\partial}{\partial\t^+}\g_\t(x,0)\,dx\\
&=&i\la \clp(\t,\la)+\g_\t(\t,0)+ \int_0^\t
e^{i\la(\t- x)}\g_\t^u(x,\t)\g_\t^l
(\t,0)\,dx\\
&=&i\la \clp(\t,\la)+\left(I_n+ \int_0^\t
e^{i\la(\t- x)}\g_\t(x,\t)\,dx\right)
\g_\t(\t,0)\\
&=&i\la \clp(\t,\la)+\left(I_n+ \int_0^\t e^{i\la
x}\g_\t(\t-x,\t)\,dx\right)
\g_\t(\t,0)\\
&=& i\la \clp(\t,\la)+\clp_*(\t,\la)\g_\t(\t,0).
\end{eqnarray*}
Here we removed the superscripts $u$ and $l$
using \eqref{paris} and using the fact that the
value of an integral does not depend on the
value of the integrand at one point. Using now
\eqref{KSeq2267} we obtain in a similar way that
\[
\frac{\partial}{\partial\t^-}\clp(\t,\la) =i\la
\clp(\t,\la)+\clp_*(\t,\la)\g_\t(\t,0).
\]
It follows that
$\frac{\partial}{\partial\t}\clp(\t,\la)$ exists
and that the first equality in
\eqref{kr1} holds.

Analogously, using the (\ref{KSeq22}) and
\eqref{KSeq2267}, we have
\begin{eqnarray*}
\frac{\partial}{\partial\t^\pm}\clp_*(\t,\la)
&=&\frac{\partial}{\partial\t^\pm}
\int_0^\t  e^{i\la x}\g_\t(\t-x,\t)\,dx\\
&=&e^{i\la\t}\g_\t(0,\t)+\int_0^\t e^{i\la
x}\frac{\partial}{\partial\t^\pm}
\g_\t(\t-x,\t)\,dx\\
&=&e^{i\la\t}\g_\t(0,\t)+\int_0^\t
e^{i\la x}\g_\t(\t-x,0)\g_\t(0,\t)\,dx\\
&=&e^{i\la\t}\left(I_r+\int_0^\t
e^{i\la (x-\t)}\g_\t(\t-x,0)\,dx\right)\g_\t(0,\t),\\
&=&e^{i\la\t}\left(I_r+\int_0^\t  e^{i\la
x}\g_\t(x,0)\,
dx\right)\g_\t(0,\t)\\
&=&\clp(\t,\la)\g_\t(0,\t).
\end{eqnarray*}

Since $k$ is hermitian, we have $\g_\t(0,\t)=\g_\t(\t,0)^*$.
Thus $\clp$ and $\clp_*$ satisfy \eqref{kr1}, and hence
$Y(\t,\la)=\begin{bmatrix}{\clp}(\t,\la)&
{\clp}_*(\t,\la)\end{bmatrix}$ satisfies \eqref{kr}. \epr

\section{Intermezzo: The continuous analogue of the resultant}
\setcounter{equation}{0}
\label{270308}
We review here the results of \cite{GKLe07b}
needed in the proof of Theorem
\ref{mainthm2}. The definition of the
resultant operator $\bfr(B,D)$ has already been given
in the introduction. Consider an entire
matrix function of the form
\begin{equation}
\label{tokyo} L(\la)=I_r+\int_0^\t e^{i\la
x}\ell(x)dx,\quad \ell\in{\mathbf L}_1^{r\times r}[0,\t].
\end{equation}
With a slight abuse of terminology, following \cite{GKLe07b}, we call $L(\la)$
a \emph{Krein orthogonal matrix function} if
there exists a hermitian ${\mathbb C}^{r\times
r}$--valued function $k\in{\mathbf L}_1^{r\times
r}[-\t, \t]$ such that
\[
\ell(t)-\int_0^\t k(t-u)\ell(u)\,du=k(t),\quad 0\le t\le
\t.
\]
In that case we refer to $\delta-k$ as
the associate \emph{weight}. The following result is proved in \cite[Theorem
5.6]{GKLe07b}.
\begin{Tm}\label{resultant1}
Let $L$ be a ${\mathbb C}^{r\times r}$-valued
entire function of the form \eqref{tokyo}. Then
there exists a hermitian matrix function
$k\in{\mathbf L}_1^{r\times r}[-\t,\t]$ such that
$L$ is the Krein orthogonal matrix function with
weight  $\delta-k$ if and only if there exists a matrix function $M$ of the form
\begin{equation}
\label{fM}
M(\la)=I_r+\int_0^\t e^{i\la u}m(u)du, \quad m\in{\mathbf L}_1^{r\times r}[0,\t],
\end{equation}
such that the following two conditions are satisfied:
\begin{eqnarray}
&& L(\la)L^\sharp(\la)=M^\sharp(\la)M(\la), \quad \la\in \BC,\label{cond1}\\
\noalign{\vskip4pt}
&& \ker L^\sharp(\la)\cap \ker
M(\la)=\left\{0\right\}, \quad \la\in \BC.
\label{cond2}
\end{eqnarray}
Furthermore, when these conditions hold, the function $k$ is
given by the formula
\begin{equation}
\label{eqweight}
k=\left[\bfr(L^\sharp,M)\right]^{-1}q,
\qquad q(u)=\begin{cases}
\ell(-u)^*,\quad -\t\le u\le0,\\
\noalign{\vskip4pt}
m(u),\quad \hspace{6mm}0\le u\le \t.
\end{cases}
\end{equation}
\end{Tm}

In \cite{GKLe07b} the above theorem  is derived as a corollary of
the following somewhat more general theorem (\cite[Theorem 5.5]{GKLe07b}).

\begin{Tm}\label{resultant2}
Given $\ell, m\in{\mathbf L}_1^{r\times r}[0,\t]$, put
\[
L(\la)=I_r+\int_0^\t e^{i\la u}\ell(u)du, \quad M(\la)=I_r+\int_0^\t e^{i\la u}m(u)du.
\]
Then there is a hermitian matrix function
$k\in{\mathbf L}_1^{r\times r}[-\t,\t]$ such that
\begin{eqnarray}
&&\ell(t)-\int_0^\t k(t-u)\ell(u)\,du=k(t),\quad 0\le t\le
\t,\label{eql}\\
\noalign{\vskip4pt}
&&m(t)-\int_0^\t m(u)k(t-u)\,du=k(t),\quad 0\le t\le
\t, \label{eqm}
\end{eqnarray}
if and only if  the  two conditions \eqref{cond1} and \eqref{cond2} are satisfied,
and in that case
the function $k$ is uniquely determined by \eqref{eqweight}.
\end{Tm}

In general, a Krein orthogonal matrix function $L$ may have many different weights.
This is reflected by the fact that given $L$ as in Theorem \ref{resultant1}
there may be many different functions $M$
of the form \eqref{fM} satisfying  \eqref{cond1} and \eqref{cond2}.
However, as soon as $M$ is fixed, then the
weight is uniquely determined by \eqref{eqweight} (as we see
from Theorem \ref{resultant2}).

\medskip
\noindent\textbf{Remark.} If in \eqref{eqweight} the functions $\ell$ and $m$ are
continuous on the interval $[0, \t]$, then the function $q$ in the
right hand side of \eqref{eqweight}
is a continuous function on $[-\t, \t]$ with a possible jump discontinuity at zero.
This implies that the function $k$ defined by  \eqref{eqweight} is also continuous
on $[-\t, \t]$ with a possible jump discontinuity at zero.

\section{Proof of Theorem \ref{mainthm2}}
\label{seckrs}
\setcounter{equation}{0}
Throughout this section $k$ is  a ${r\ts r}$-matrix valued
accelerant on $[-\mathbf
T,\mathbf T]$, with possibly a jump
discontinuity at the origin, and we consider the Krein system \eqref{kr}
with the potential $a$ defined by $k$. Furthermore $U(\t,\la)$
will be the matrizant of \eqref{kr}.

Our aim is to prove Theorem \ref{mainthm2}. As in Theorem \ref{mainthm2}, put
\[
F(\la)=e^{i\la\bft}\begin{bmatrix} I_r&I_r
\end{bmatrix}U(\bft, -\la)\begin{bmatrix}
I_r\\0
\end{bmatrix}, \  G(\la)=\begin{bmatrix} I_r&I_r
\end{bmatrix}U(\bft, -\la)\begin{bmatrix}
0\\I_r
\end{bmatrix}.
\]
First let us show that
\begin{equation}
\label{basicids}
F(\la)=e^{i\la\bft}\clp(\bft, -\la) \quad \mbox{and}\quad G(\la)=\clp_* (\bft, -\la),
\end{equation}
where $\clp(\bft, \la)$ and $\clp_*(\bft, \la)$ are defined by \eqref{defP_mat}
and \eqref{defP*_mat} with $\t=\bft$. To obtain \eqref{basicids}
note that for each $\la\in \BC$ the two $r \ts 2r$
matrix functions
\[
\begin{bmatrix} I_r &I_r\end{bmatrix}U(\t,\la) \quad \mbox{and} \quad \begin{bmatrix}{\mathcal
P}(\t,\la)&
{\mathcal P}_*(\t,\la)\end{bmatrix}
\]
satisfy the linear differential equation
\eqref{kr}, and  at $\t=0$ both functions are equal to $\begin{bmatrix} I_r &I_r\end{bmatrix}$. Thus both
have the same initial
condition at $\t=0$.
It follows that these two functions coincide on $0\leq \t\leq \bft$. For $\t=\bft$ this yields
the identities in \eqref{basicids}.

Using \eqref{basicids}, we see from the formulas for $\clp$ and $\clp_*$ in
\eqref{defP_mat} and \eqref{defP*_mat} that
\begin{equation}
\label{eqFG}
F(\la)= I_r+\int_0^{\bft}f(x)e^{i\la
x}dx, \quad G(\la)=I_r+\int_{-\bft}^0 g(x)e^{i\la
x}dx,
\end{equation}
with $f(x)=\g_\bft(x,0)$ on $0\leq x\leq \bft$ and
$g(x)=\g_\bft(\bft+x,\bft)$ on the interval $-\bft\leq x\leq 0$. In particular, the
functions $f$ and $g$ are continuous on their respective domains as desired.

 It remains to prove \eqref{Pernety}. To do this we first
 derive the following lemma.

\begin{La}\label{lemfact}
The functions $\clp$ and
$\clp_*$ given by \eqref{defP_mat} and
\eqref{defP*_mat}, respectively, satisfy the
identity
\begin{equation}
\label{rlcanonfact}
\clp(\t, \la)\clp^\sharp(\t, \la)=\clp_*(\t,
\la)\clp_*^\sharp(\t, \la)\quad (0\leq \t\leq \bft,\
\la\in \BC).
\end{equation}
Furthermore,  for each $0\leq \t\leq \bft$ the
left hand side in the above identity is a right
canonical factorization (that is, for each $0\leq \t\leq \bft$ the function
$\det \clp(\t, \la)$ has no zero in the closed lower half plane) while the right side is
a left canonical factorization (that is, for each $0\leq \t\leq \bft$ the function
$\det \clp_*(\t, \la)$ has no zero in the closed upper  half plane).
\end{La}
\bpr Fix $0\leq \t\leq \bft$. Recall that  the integral operator
$T_\t$ defined by \eqref{intor} is selfadjoint and invertible.
Let $a_\t$, $b_\t$, $b_\t$, $d_\t$
be the $L^1$-functions defined by
\begin{eqnarray*}
&&a_\t(t)- \int_0^\t k(t-u)a_\t(u)\, du=k(t),
\quad 0\leq t\leq \t,\\
&&b_\t(t)- \int_{-\t}^0 b_\t(u)k(t-u)\,
du=k(t), \quad -\t\leq t \leq 0,\\
&& c_\t(t)-\int_{-\t}^0 k(t-u)c_\t(u)\, du=k(t),
\quad -\t\leq t \leq 0,\\
&&d_\t(t)-\int_0^\t d_\t(u)k(t-u)\, du=k(t),
\quad 0\leq t\leq \t,
\end{eqnarray*}
and put
\begin{eqnarray*}
&&\cla_\t(\la)= I+\int_0^\t e^{i\la
s}a_\t(s)\,ds,\quad \clb_\t(\la)= I+\int_{-\t}^0
e^{i\la s}b_\t(s)\,ds
,\\
&&\clc_\t(\la)= I+\int_{-\t}^0 e^{i\la
s}c_\t(s)\,ds,\quad \cld_\t(\la)= I+\int_0^\t
e^{i\la s}d_\t(s)\,ds.
\end{eqnarray*}
In terms of the resolvent kernel $\g_\t(t,s)$ associated with $k$ we have
\begin{eqnarray*}
&&a_\t(x)=\g_\t(x, 0),\hspace{.9cm} b_\t(-x)=\g_\t(0,x)\quad (0\leq x\leq \t);\\
&&c_\t(x)=\g_\t(\t+x, \t)\quad d_\t(-x)=\g_\t(\t,\t+x)\quad (-\t\leq x\leq 0).
\end{eqnarray*}
Note that in this terminology, the functions
$\clp$ and $\clp_*$ given by (\ref{defP_mat}) and
(\ref{defP*_mat}) are equal to
\begin{equation}
\label{PP*} \clp(\t,\la)=e^{i\la\t}\cla_\t(-\la),
\quad \clp_*(\t,\la)=\clc_\t(-\la).
\end{equation}
From Theorem 5.3 in \cite{GKLe07b} we know that
\begin{equation}
\label{qcpker}
\cla_\t(\la)\clb_\t(\la)=\clc_\t(\la)\cld_\t(\la),
\quad \kr\clb_\t(\la)\cap\kr\cld_\t(\la)=\{0\}.
\end{equation}
Next recall that $k$ is hermitian. This implies
that
\[
b_\t(-x)=a_\t(x)^*,\quad c_\t(-x)=d_\t(x)^*
\quad (0\leq x\leq \t),
\]
and hence $\cla_\t^\sharp(\la)=\clb_\t(\la)$ and
$\cld_\t^\sharp(\la)=\clc_\t(\la)$. In particular,
(\ref{qcpker}) reduces to
\begin{equation}
\label{qcpker*}
\cla_\t(\la)\cla_\t^\sharp(\la)=\cld_\t^\sharp(\la)\cld_\t(\la),
\quad \kr\cla_\t^\sharp(\la)\cap\kr\cld_\t(\la)=\{0\}.
\end{equation}

Finally, since for each $0\leq \t \leq
\bft$ the operator $T_\t$ in (\ref{intor}) is
selfadjoint and invertible, it follows that $T_\t$ is
strictly positive for each $0\leq \t \leq
\bft$. Then we know (using the theory
of Krein orthogonal functions; see Theorem 8.1.1
in \cite{ellis-gohberg}) that
the function $\det\cla_\t(\la)$ has no zero in the closed upper
half plane, and the function $\det\cld_\t^\sharp(\la)$ has no zero in
the closed lower half plane. Thus
$\cla_\t(\la)\cla_\t^\sharp(\la)$ is a  left
canonical factorization and
$\cld_\t^\sharp(\la)\cld_\t(\la)$ is a right canonical factorization.
Using (\ref{PP*}) the above remarks provide the
proof of the lemma. \epr

\medskip
We are now ready to prove \eqref{Pernety}. From \eqref{basicids}
and \eqref{rlcanonfact} it follows that
\begin{equation}
\label{GFlr}
F(\la)F^\sharp(\la)=G(\la)G^\sharp(\la).
\end{equation}
Moreover the left hand side of this identity is a left canonical factorization and the right
hand side is a right canonical factorization. In particular, $\kr F^\sharp\cap\kr G^\sharp=\{0\}$.
This allows us to apply Theorem \ref{resultant2} with $\t=\bft$, $\ell(u)=f(u)$ and $m(u)=g(-u)^*$,
where the functions $f$ and $g$ are as in (\ref{eqFG}). In other words, we apply
$L=F$
and $M=G^\sharp$. It follows that there exists a unique hermitian $\tilde{k}\in{\mathbf L}_1^{r\times
r}[-\bft, \bft]$ such that
\begin{eqnarray}
&&f(t)-\int_0^\bft \tilde{k}(t-s)f(s)\,ds=\tilde{k}(t),\quad 0\leq t\leq
\bft,\label{eqf}\\
\noalign{\vskip4pt}
&&g(t)-\int_{-\bft}^0 \tilde{k}(t-s)g(s)\,ds=\tilde{k}(t),\quad -\bft\le t\leq 0.\label{eqg}
\end{eqnarray}
Moreover, $\tilde{k}$ is given by the formula
\begin{equation*}
\tilde{k}=[\bfr(F^\sharp,G^\sharp)]^{-1}q  \quad\mbox{with} \quad q(x)=\begin{cases}
f(-x)^*,\quad -\bft\le x\le0,\\
g(-x)^*,\quad \hspace{4mm}0\le x\le \bft.
\end{cases}
\end{equation*}
Since $f(x)=\g_\bft(x,0)$ on $0\leq x\leq \bft$ and
$g(x)=\g_\bft(\bft+x,\bft)$ on the interval $-\bft\leq x\leq 0$, we know from
the proof of Lemma \ref{lemfact} that \eqref{eqf} and \eqref{eqg} also hold with $\tilde{k}$
being replaced by the original accelerant $k$. But then, by the uniqueness statement
in Theorem \ref{resultant2}, the functions  $\tilde{k}$ and $k$ coincide. Thus \eqref{Pernety}
holds, which completes the proof of Theorem \ref{mainthm2}.

\medskip\noindent
\textbf{Remark.} In the proof of Lemma \ref{lemfact} we used in an essential way the accelerant and its
properties. However, this is not necessary. It is possible to give a proof of Lemma 4.1 without
any reference to the accelerant. In fact, such a proof can be obtained by using the properties of
a canonical differential systems of Dirac type. To see this note
that $e^{-i\t \la}Y(\t, -2\bar{\la})^*$ is a solution of a  canonical differential system of
Dirac type with potential $v(\t)=-ia(\t)$  whenever $Y(\t, \la)$ is a solution of \eqref{kr}.
We will come back to this in a later paper.

\section{An example with jump discontinuity: the
rational case}\label{exjump}
\setcounter{equation}{0}
\label{Felix_Faure}
In this section we consider the case where the
accelerant is of the form
\begin{equation}\label{place_de_la_madeleine}
k(t)=\begin{cases}iCe^{-itA}
(I-P)B,\quad t>0,\\
-iCe^{-itA}PB,\,\,\,\,\qquad t<0.
\end{cases}
\end{equation}
In this expression, $A,B$ and $C$ are matrices
of appropriate sizes and $P$ is a projection
commuting with $A$. Motivation for such a form
originates with linear system theory. Indeed, let
$W$ be a rational ${\BC}^{p\times
q}$-valued function, analytic at infinity. Then,
as is well-known, $W$ admits a realization of
the form
\[
W(\la)=D+C(\la I_N-A)^{-1}B,\]
where $D=W(\infty)$ and $(A,B,C)\in{\mathbb
C}^{N\times N}\times{\mathbb C}^{N\times q}\times
{\mathbb C}^{p\times N}$.
Assume furthermore that $A$ has no real
eigenvalues. Then, the function $W$ belongs to
the Wiener algebra, and
\[
W(\la)=D+\int_{\BR}e^{i\la t}k(t)dt,
\]
where $k$ is of the form
\eqref{place_de_la_madeleine} with $P$ being  the
Riesz projection corresponding to the eigenvalues
of $A$ in the upper--half plane. Note that, in
general, functions $k$ of the form
\eqref{place_de_la_madeleine} need not have
summable entries.

\medskip
In this section we first take
\begin{equation}
\label{matABC}
A=\begin{bmatrix}a^\times &-bb^*\\0&a^{\times
*}\end{bmatrix},\quad B=\begin{bmatrix}b\\
c^*\end{bmatrix},\quad
C=\begin{bmatrix}-c&-b^*\end{bmatrix},
\end{equation}
where $(a,b,c)\in{\mathbb C}^{n\times n}\times
{\mathbb C}^{n\times k}\times{\mathbb
C}^{k\times n}$, and throughout it is assumed that the spectra
of $a$ and $a^\times =a-bc$ are both in the open
upper half-plane. For $P$ we take the Riesz projection of $A$
corresponding to the eigenvalues
in the upper--half plane. In other words  $P$ is given by
\begin{equation}
\label{matP}
P=\begin{bmatrix}
I&i\om\\
0&0
\end{bmatrix},
\end{equation}
where $\om$ is the unique solution of the Lyapunov equation
\begin{equation}
\label{defOm}
i(\Omega a^{\times *}-a^\times \Omega)=bb^*.
\end{equation}
With $A, B, C$ and $P$ as in \eqref{matABC} and \eqref{matP}, the function
\begin{equation}\label{symb}
W(\la)=I_r+\int_{\BR}e^{i\la t}k(t)dt
\end{equation}
is
positive definite on the real line. Conversely, any rational $r\ts r$ matrix function $W$
which is positive definite on the real line and analytic at infinity with $W(\iy)=I_r$
can be represented in this way (see \cite{ag3}).

\begin{Pn}
\label{strictpsexp}When $k$ is of the form
\eqref{place_de_la_madeleine}  with $A, B, C$ and $P$ being given by \eqref{matABC} and \eqref{matP},
then $k$ is an accelerant on each interval $[-\bft,\bft]$. Moreover, in this case the
corresponding potential is given by
\begin{equation}
\label{formpot}
a(\t)=i\left( (I_n+\Omega(Y-e^{-i\t a^*}Ya^{i\t
a}))^{-1}(b+i \Omega c^*)\right)^*,
\end{equation}
where $\Omega $ is given by \eqref{defOm}, and where $Y$ is the  solution of the
Lyapunov equation
\begin{equation}
\label{defY}
i(Ya-a^*Y)=-c^*c.
\end{equation}
\end{Pn}
\bpr The fact that the function
$W$ in \eqref{symb} is positive definite on the real line implies that for each $\t$
the integral operator $T_\t$ in \eqref{intor} is strictly positive. Hence $k$
is an accelerant on each interval $[-\bft,\bft]$.

Using Theorem 4.1 in \cite{bgk3} one computes that in this setting
\[
\gamma_\t(0,\t)=-iC(Pe^{-i\t (A-BC)}\big|_{\im P})^{-1}PB.
\]
Since $A, B$ and $C$ are given by \eqref{matABC}, we have
\[
A-BC=\begin{bmatrix}
a&0\\
c^*c&a^*
\end{bmatrix}.
\]
It then follows,
as computed in \cite[p.15]{ag3}, that
\[
\gamma_\t(0,\t)=-i(I_n+\Omega(Y-e^{-i\t
a^*}Ya^{i\t a}))^{-1}(b+i\Omega c^*),
\]
where $\Omega $ and $Y$ are given by \eqref{defOm} and \eqref{defY}, respectively.
Since the potential is given by $a(\t)=\g_\t(\t,0)$ and $\g_\t(\t,0)=\g_\t(0, \t)^*$,
we see that $a$ is given by \eqref{formpot}. \epr

\medskip

Next we assume that the matrices $A$, $B$, and $C$ in \eqref{place_de_la_madeleine} are given by
\begin{equation}\label{ABC2}
A=2\begin{bmatrix}
\b^*&\g_2\g_2^*\\
0&\b
\end{bmatrix},\quad B=\sqrt{2}\begin{bmatrix}
\g_2\\
\g_1
\end{bmatrix},\quad C=\sqrt{2}\begin{bmatrix}
\g_1^*&\g_2^*
\end{bmatrix},
\end{equation}
where $\b$ is a square matrix of order $n$, and $\g_1$ and $\g_2$
are matrices of sizes $n\ts r$. Furthermore, we assume that $\b^*-\b=i\g_2\g_2^*$.
A triple of matrices $\b$, $\g_1$ and $\g_2$ with these properties will be called \emph{admissible}.
For the matrix
$P$ in \eqref{place_de_la_madeleine} we take
\begin{equation}\label{matP2}
P=\begin{bmatrix}
I_n&-i I_n\\
0&0
\end{bmatrix}.
\end{equation}
The fact that the triple of matrices $\b$, $\g_1$ and $\g_2$ is assumed to be admissible
implies that with $A, B, C$ and $P$ as in \eqref{ABC2} and \eqref{matP2}, the function
\eqref{symb} is positive semi-definite on the real line. Conversely, any rational $r\ts r$ matrix function $W$
which is positive semi-definite on the real line and analytic at infinity with $W(\iy)=I_r$
can be represented in this way (see \cite{GKSakh98}, also \cite{GKSakh02}).

For information about the connection
between the matrices $A$, $B$, $C$ and $P$ in \eqref{matABC} and \eqref{matP} and those
in \eqref{ABC2} and \eqref{matP2}, we refer to the introduction of
\cite{agks00}.

\begin{Pn}
\label{psexpo} Let $\b$, $\g_1$ and $\g_2$ be an admissible triple, and put
\begin{equation}
\label{kernelpsexp} k(t)=-2(\g_1+i\g_2)^*e^{-2it\b}\g_1, \quad k(-t)=k(t)^* \quad (t>0).
\end{equation}
Then $k$ is an accelerant on each interval $[-\bft,\bft]$, and  the
corresponding potential is given by
\begin{equation}
\label{potpsexp} a(\t)=-2(\g_1+i\g_2)^*e^{-i\t
\a^*}\Sigma(\t)^{-1}e^{-i\t \a}\g_1, \quad \a=\b-\g_1\g_2^*,
\end{equation}
where
\begin{equation}
\label{SigLam} \Sigma(t)=I_n+\int_0^t\Lambda(s)\Lambda(s)^*\, ds, \quad \Lambda(t)=\begin{bmatrix}
e^{-it\a}\g_1&-e^{it\a}(\g_1+i\g_2)
\end{bmatrix}.
\end{equation}
\end{Pn}
\bpr
Let $A, B, C$ and $P$ be given by  \eqref{ABC2} and \eqref{matP2}. Put
\[
S= \begin{bmatrix}
I_n&i I_n\\
0&I_n
\end{bmatrix}.
\]
Then $S$ is invertible, and one computes that
\begin{eqnarray*}
&&A =S\begin{bmatrix}
2\b^*&0\\
0&2\b
\end{bmatrix}S^{-1}, \quad B=\sqrt{2}S\begin{bmatrix}
-i(\g_1+i\g_2)\\
\g_1
\end{bmatrix},\\
\noalign{\vskip6pt}
&&C=\sqrt{2}\begin{bmatrix}\g_1^*&i(\g_1+i\g_2)^*\end{bmatrix}S^{-1},\quad
P= S\begin{bmatrix}
I_n&0\\
0&0
\end{bmatrix}S^{-1}.
\end{eqnarray*}
It follows that
\begin{eqnarray*}
iCe^{-itA}(I-P)B&=& 2i\begin{bmatrix}\g_1^*&i(\g_1+i\g_2)^*\end{bmatrix}\begin{bmatrix}
0&0\\
0&e^{-2it\b}
\end{bmatrix}\begin{bmatrix}
-i(\g_1+i\g_2)\\
\g_1
\end{bmatrix}\\
&=& -2(\g_1+i\g_2)^*e^{-2it\b}\g_1.
\end{eqnarray*}
Analogously
\begin{eqnarray*}
-iCe^{-itA}PB&=& -2i\begin{bmatrix}\g_1^*&i(\g_1+i\g_2)^*\end{bmatrix}\begin{bmatrix}
e^{-2it\b^*}&0\\
0&0
\end{bmatrix}\begin{bmatrix}
-i(\g_1+i\g_2)\\
\g_1
\end{bmatrix}\\
&=& -2\g_1^*e^{-2it\b^*}(\g_1+i\g_2).
\end{eqnarray*}
It follows that $k$ given by \eqref{kernelpsexp} can be written in the form \eqref{place_de_la_madeleine}
with $A, B, C$, and $P$ as in \eqref{ABC2} and \eqref{matP2}.

Next, we consider $A^\ts=A-BC$. We have
\[
A^\ts = 2\begin{bmatrix}
\a^*&0\\
-\g_1\g_1^*&\a
\end{bmatrix}, \quad\mbox{where $\a=\b-\g_1\g_2^*$}.
\]
The proof of Proposition 4.1 in \cite{GKSakh98} shows that
\[
Pe^{-itA^\ts}|_{\im P}=e^{-it\a}\Sigma(t)e^{-it\a^*}, \quad t\geq 0.
\]
Since $\Sigma(t)$ is positive definite, the matrix $\Sigma(t)$ is
invertible. Hence the map $Pe^{-itA^\ts}|_{\im P}$, viewed as an
operator acting on $\im P$, is invertible. By Theorem 4.3 in
\cite{bgk3} this implies that for our $k$ the integral operator
$T_\t$ given  by \eqref{intor} is invertible for each $\t$, and
\begin{eqnarray*}
\gamma_\t(0,\t)&=&-iC(Pe^{-i\t (A-BC)}\big|_{\im P})^{-1}PB\\
&=& -2 \g_1^*e^{i\t\a^*}\Sigma(\t)^{-1}e^{i\t\a}(\g_1+i\g_2).
\end{eqnarray*}
Here we used that
\[
PB= \sqrt{2}\begin{bmatrix}
-i(\g_1+i\g_2)\\0
\end{bmatrix}, \quad C|_{\im P}= \sqrt{2}\g_1^*.
\]
Since the potential is given by $a(\t)=\g_\t(\t,0)$ and $\g_\t(\t,0)=\g_\t(0, \t)^*$,
we see that $a$ is given by \eqref{potpsexp}. \epr

\medskip
\noindent\textbf{Remark.}
Note that the two propositions in this section do not cover the example
presented in the introduction. Indeed, when $k$ is given by \eqref{Ranelagh},
then $k$ is not an accelerant for $[-\pi/2, \pi/2]$.

\section{Another class of potentials}\label{otherclass}
\setcounter{equation}{0}

We now consider the case where the $\BC^{r\ts r}$-valued  accelerant $k$
admits a representation of the form
\begin{equation}
\label{voltaire} k(t)=Ce^{tA}B,\quad
t\in[-{\mathbf T},{\mathbf T}],
\end{equation}
where $A,B$ and $C$ are matrices of appropriate sizes. We assume that
there exists a  hermitian matrix $H$ such that
\begin{equation}\label{condH}
HA+A^*H=0 \quad \mbox{and} \quad C=B^*H.
\end{equation}
The latter implies that $k(t)^*=k(-t)$ on $[-{\mathbf T},{\mathbf T}]$, and hence
$k$ is a hermitian kernel. Under certain minimality
conditions the converse statement is also true. More precisely, if $k$ given by
\eqref{voltaire} with the pair $(A,B)$ being  controllable and the pair $(C,A)$ being
observable, then $k(t)^*=k(-t)$ implies that there exists a unique invertible
hermitian matrix $H$ such that \eqref{condH}  holds.

Let $\tau\in(0,{\mathbf T}]$. As proved in
\cite{GKvS04}, equation \eqref{langoustine} has a
unique solution if and only if the matrix
\begin{equation}
\label{Mtau}
M_{\tau}=I-\int_0^\tau e^{-sA}BCe^{sA}ds
\end{equation}
is invertible. When this is the case, we have:
\begin{equation}
\label{montparnasse}
\gamma_{\tau}(t,s)=Ce^{tA}M_{\tau}^{-1}e^{-sA}B,
\quad s,t\in[-{\tau},{\tau}].
\end{equation}

\begin{Pn}\label{expkernel} Assume $k$ is given by \eqref{voltaire},
and let $H$ be a hermitian matrix $H$ such that \eqref{condH}  holds.
Then $k$ is an accelerant if and only if  the matrix $M_{\tau}$ in \eqref{Mtau}
is non-singular for $0\leq \t\leq \bft$. In that case
 the corresponding potential  is given by
\begin{equation}\label{potenexp}
a(t)=Ce^{tA}M_{t}^{-1}B,\quad 0<t\leq \bft,
\end{equation}
and the functions
\begin{eqnarray*}
\clp(\t, \la)&=&e^{i\la\t}I_r+C(A-i\la I_r)^{-1}\left\{e^{\tau
  A}-e^{i\la\t}I_r\right\}M_\t^{-1}B,\\
  \noalign{\vskip4pt}
\clp_*(\t, \la)&=&I_r+C(A-i\la I_r)^{-1}\left\{e^{\tau
  A}-e^{i\la\t}I_r\right\}M_\t^{-1}e^{-\t A}B,
  \end{eqnarray*}
are the associate  Krein orthogonal matrix functions.
\end{Pn}
\bpr Since $k$ is hermitian, the operator $T_\t$ will be strictly
positive if and only $T_\t$ is invertible. The latter happens if and
only if $M_\t$ is non singular. Thus  $k$ is an accelerant if and
only if $M_\t$ is non-singular for $0\leq \t\leq \bft$.

Assume $k$ to be  an accelerant. Then the potential is given by
$a(t)=\gamma_{t}(t,0)$ on $(0,\bft]$.
Using \eqref{montparnasse}, this yields \eqref{potenexp}. Furthermore,
the associate Krein orthogonal function $\clp$ for $k$ can be
computed as follows:
\begin{eqnarray*}
\clp(\t,\la)&=& e^{i\la\t}\left(I_r+\int_0^\t e^{-i\la x}\gamma_\t(x,0)
dx\right)\\
\noalign{\vskip4pt}
&=&e^{i\la\t}\left(I_r+\int_0^\t e^{-i\la x}Ce^{xA}M_\t^{-1}Bdx\right)\\
\noalign{\vskip4pt}
&=&e^{i\la\t}\left(I_r+C\left(\int_0^\t
e^{-i\la
  x}e^{xA}dx\right)M_\t^{-1}B\right)\\
\noalign{\vskip4pt}
&=&e^{i\la\t}\left(I_r+C(A-i\la
I)^{-1}\left\{e^{\tau
(A-i\la)}-I_r\right\}M_\t^{-1}B\right)\\
\noalign{\vskip4pt}
&=&e^{i\la\t}I_r+C(A-\la I)^{-1}\left\{e^{\tau
A}-e^{i\la\t}I_r\right\}M_\t^{-1}B.
\end{eqnarray*}
Analogously,
\begin{eqnarray*}
\clp_*(\t,\la)&=&I_r+\int_0^\t e^{i\la x}\gamma_\t(\t-x,\t)
dx\\
\noalign{\vskip4pt}
&=&I_r+\int_0^\t e^{i\la x}Ce^{(\t-x)A}M_\t^{-1}e^{-\t A}Bdx\\
\noalign{\vskip4pt}
&=&I_r+Ce^{\t A}\left(\int_0^\t e^{i\la x} e^{-xA}dx\right)M_\t^{-1}e^{-\t A}B\\
\noalign{\vskip4pt}
&=&I_r+Ce^{\t A}\left((i\la -A)^{-1}e^{(i\la -A)\t}-(i\la -A)^{-1}\right)M_\t^{-1}e^{-\t A}B\\
\noalign{\vskip4pt}
&=&I_r+C(A-i\la I_r)^{-1}\left(e^{\t A}-e^{i\la \t}I_r \right)M_\t^{-1}e^{-\t A}B.
\end{eqnarray*}
This completes the proof. \epr

\begin{Cy}
Assume $k$ is given by \eqref{voltaire},
and assume that \eqref{condH}  holds with $H=-I$. Then $k$ is an accelerant.
In particular, if  $r_j>0$ and $\b_j\in \BR$ for $j=1, \ldots, n$, then  the function
\begin{equation}
\label{exk} k(t)=-\sum _{\nu=1}^n r_n e^{i\b_\nu t}
\end{equation}
is an accelerant for each each interval $[-{\bft},{\bft}]$.
\end{Cy}
\bpr From $H=-I$, we see that the matrix $M_\t$ in \eqref{Mtau} can be rewritten as
\[
M_\t=I+\int_0^\tau \left(Ce^{sA}\right)^*\left(Ce^{sA}\right)ds.
\]
It follows that $M_\t$ is positive definite and hence non-singular for each $\t\geq 0$. Thus $k$ is an
accelerant by Proposition \ref{potenexp} above.

Next, consider the function $k$ in \eqref{exk}. Since $r_j>0$ and $\b_j\in \BR$ for each $j=1, \ldots, n$,
we can represent $k$ as in \eqref{voltaire} by taking
\[
A=\diag(i\b_1, i\b_2, \ldots, i\b_n), \quad
C=\begin{bmatrix}\sqrt{r_1} &\sqrt{r_2}& \cdots &\sqrt{r_n} \end{bmatrix}, \quad B=-C^*.
\]
But then \eqref{condH} holds with $H=I$. By the result of the first paragraph,
this shows that $k$ is an
accelerant on $[-{\bft},{\bft}]$ for each $\bft>0$.\epr

\medskip
From \eqref{Mtau} it follows that
\begin{equation}
\label{Mprime}M_\tau^\prime=\frac{d}{d\t}M_\tau=-e^{-\t A}BCe^{\t A}.
\end{equation}
This together with the explicit formula \eqref{montparnasse} allows
one to give a direct proof of the Krein-Sobolev
equation \eqref{KSeq1} and of equation
\eqref{KSeq2} for accelerants as in \eqref{voltaire}.

The class of accelerants considered in this section includes the
restrictions of polynomials to $[-{\mathbf T},{\mathbf T}]$. On  the
other hand, when considered for $t$ on the whole real line, $k$ is
never integrable (except for the trivial case  $k=0$). Thus this
class of accelerants has a zero intersection with the accelerants
considered in the first part of the previous section. Nevertheless
the class of potentials corresponding to the accelerants considered
in this section shares a number of common properties with the
strictly pseudo-exponential potentials. For instance, using
\eqref{Mprime}, we have
\[
\begin{split}
a(0)&=CB\\
a^\prime(0)&=CAB+(CB)^2\\
&\hspace{2mm}\vdots
\end{split}
\]
and there exist non commutative polynomials
$f_0, f_1,\ldots$ such that
\[
CA^\ell B=f_\ell(v(0),\ldots,
v^{(\ell)}(0)),\quad \ell=0,1,\ldots
\]
Thus, and as for strictly pseudo-exponential potentials
(see \cite{ag9}), one can in principle recover
the potential from the values of its first
derivatives at the origin (cf., \cite{GKSakh05}, where such results are proved for
pseudo-exponential potentials).

\end{document}